\documentclass[11pt]{article}

\usepackage{geometry,graphicx,latexsym,amssymb,verbatim, multicol}
\usepackage{tikz}
\usetikzlibrary{calc}

\newtheorem{thm}{Theorem}

\newtheorem{prob}{Problem}
\newtheorem{cor}[thm]{Corollary}
\newtheorem{ob}[thm]{Observation}

\newtheorem{conj}{Conjecture}

\newcommand{\dontshow}[1]{}

\newcommand{\sd}{\gamma_{\rm sd}}

\newcommand{\proof}{\noindent\textbf{Proof. }}

\newcommand{\qed}{$\Box$}

\newcommand{\mod}{{\rm mod}}

\newcommand{\NH}{{\rm NH}}

\newcommand{\de}{\delta}

\newcommand{\bard}{{\overline{d}}}

\newcommand{\TR}[1]{\mbox{$\tau(#1)$}}

\newcommand{\TRt}[1]{\mbox{$\tau_t(#1)$}}

\newcommand{\3}{ \vspace{0.3cm} }
\newcommand{\2}{ \vspace{0.2cm} }

\newcommand{\cG}{{\cal G}}

\newcommand{\oo}{{\rm o}}

\def\vertex(#1){\put(#1){\circle*{2}}}
\def\vertexo(#1){\put(#1){\circle{2}}}
\def\vert(#1){\put(#1){\circle*{1.5}}}
\def\verto(#1){\put(#1){\circle{1.5}}}
\def\lab(#1)#2{\put(#1){\makebox(0,0)[c]{#2}}}
\begin{document}

\title{Simultaneous Domination in Graphs}

\author{$^1$Yair Caro and $^2$Michael A. Henning\footnote{Research supported in part by the University of Johannesburg and the South African National Research Foundation.} \\
\\
$^1$Department of Mathematics\\
University of Haifa-Oranim \\
Tivon 36006, Israel \\
Email: yacaro@kvgeva.org.il \\
\\
$^2$Department of Mathematics \\
University of Johannesburg \\
Auckland Park 2006, South Africa \\
Email: mahenning@uj.ac.za }

\date{}
\maketitle

\begin{abstract}
Let $F_1, F_2, \ldots, F_k$ be graphs with the same vertex set $V$. A
subset $S \subseteq V$ is a simultaneous dominating set if for every
$i$, $1 \le i \le k$, every vertex of $F_i$ not in $S$ is adjacent to
a vertex in $S$ in $F_i$; that is, the set $S$ is simultaneously a
dominating set in each graph $F_i$. The cardinality of a smallest
such set is the simultaneous domination number.
%We investigate bounds on the simultaneous domination number.
We present general upper bounds on the simultaneous domination
number. We investigate bounds in special cases, including the cases
when the factors, $F_i$, are $r$-regular or the disjoint union of
copies of $K_r$. Further we study the case when each factor is a
cycle.
\end{abstract}


%\newpage
\section{Introduction}

Given a collection of graphs $F_1, \ldots, F_k$ on the same vertex
set~$V$, we consider a set of vertices which dominates all the graphs
simultaneously. This was first explored by Brigham and
Dutton~\cite{bd} who defined such a set as a factor dominating set
and by Sampathkumar~\cite{Sam} who used the name global dominating
set. The natural question is what is the minimum size of a
simultaneous dominating set. This question has been studied
in~\cite{bc1,CaYu00,DaGoHeLa06,DaLa03} and \cite[Section 7.6]{hhs1}
and elsewhere. In this paper we will use the term ``simultaneous
domination" rather than ``global domination" (see~\cite{bc1,Sam}) or
``factor domination" (see~\cite{bd,DaGoHeLa06,DaLa03}).

A \emph{dominating set} of $G$ is a set $S$ of vertices of $G$ such
that every vertex outside $S$ is adjacent to some vertex in $S$.
%; that is, $N[S] = V$.
The \emph{domination number} $\gamma(G)$ is the minimum cardinality
of a dominating set in $G$. For $k \ge 1$, a $k$-\emph{dominating
set} of $G$ is a set $S$ of vertices of $G$ such that every vertex
outside $S$ is adjacent to at least~$k$ vertices in $S$. For a survey
see~\cite{hhs1,hhs2}.

Following the notation in~\cite{DaGoHeLa06}, we define a
\emph{factoring} to be a collection $F_1, F_2,\ldots, F_k$ of (not
necessarily edge-disjoint) graphs with common vertex set $V$ (the
union of whose edge sets is not necessary the complete graph). The
\emph{combined graph} of the factoring, denoted by
$G(F_1,\ldots,F_k)$, has vertex set $V$ and edge set $\bigcup_{i=1}^k
E(F_i)$. We call each $F_i$ a \emph{factor} of the combined graph.

A subset $S \subseteq V$ is a \emph{simultaneous dominating set},
abbreviated SD-set, of $G(F_1,\ldots,F_k)$ if $S$ is simultaneously a
dominating set in each factor $F_i$ for all $1 \le i \le k$. We
remark that in the literature a SD-set is also termed a \emph{factor
dominating set} or a \emph{global dominating set}. The minimum
cardinality of a SD-set in $G(F_1,\ldots,F_k)$ is the
\emph{simultaneous domination number} of $G(F_1,\ldots,F_k)$, denoted
by $\sd(F_1,F_2,\ldots,F_k)$.
%For $r \ge 1$, a subset $S \subseteq V$ is a \emph{simultaneous
%$r$-dominating set}, abbreviated SrD-set, of $G(F_1,\ldots,F_k)$ if
%$S$ is simultaneously an $r$-dominating set in each factor $F_i$ for
%all $1 \le i \le k$. The minimum cardinality of an SrD-set in
%$G(F_1,\ldots,F_k)$ is the \emph{simultaneous $r$-domination number}
%of $G(F_1,\ldots,F_k)$, denoted by $\sd^r(F_1,F_2,\ldots,F_k)$. In
%particular, we note that when $r = 1$, an SrD-set is a SD-set and
%$\sd^r(F_1,F_2,\ldots,F_k) = \sd(F_1,F_2,\ldots,F_k)$.
%
We remark that the notion of simultaneous domination is closely related to the notion of colored domination studied, for example, in~\cite{Color_Dom} and elsewhere.

For $k \ge 2$ and $\delta \ge 1$, let $\cG_{k,\delta,n}$ be the
family of all combined graphs on $n$ vertices consisting of $k$
factors each of which has minimum degree at least~$\delta$ and define
\[
\sd(k,\delta,n) = \max \{ \sd(G) \mid G \in \cG_{k,\delta,n} \}
\]

%Further let $1 \le r \le \delta$ and define
%\[
%\sd(k,\delta,n) = \max \{ \sd(G) \mid G \in \cG_{k,\delta,n} \} \hspace*{0.5cm} \mbox{and} \hspace*{0.5cm}
%\sd^r(k,\delta,n) = \max \{ \sd^r(G) \mid G \in \cG_{k,\delta,n} \}.
%\]
%
For notational convenience, we simply write $\sd(k,n) = \sd(k,1,n)$.

\subsection{Graph Theory Notation and Terminology}

For notation and graph theory terminology, we in general
follow~\cite{hhs1}. Specifically, let $G$ be a graph with vertex set
$V(G)$ of order~$n = |V(G)|$ and edge set $E(G)$ of size~$m =
|E(G)|$. The \emph{open neighborhood} of a vertex $v \in V(G)$ is
$N_G(v) = \{u \in V(G) \, | \, uv \in E(G)\}$ and the \emph{closed
neighborhood of $v$} is $N_G[v] = N_G(v) \cup \{v\}$. For a set $S
\subseteq V(G)$, its \emph{open neighborhood} is the set $N(S) =
\bigcup_{v \in S} N(v)$ and its \emph{closed neighborhood} is the set
$N[S] = N(S) \cup S$. The degree of $v$ is $d_G(v) = |N_G(v)|$. Let
$\delta(G)$, $\Delta(G)$ and $\bard(G)$ denote, respectively, the
minimum degree, the maximum degree and the average degree in $G$. If
$d_G(v) = k$ for every vertex $v \in V$, we say that $G$ is a
$k$-regular graph. If the graph $G$ is clear from the context, we
simply write $N(v)$, $N[v]$, $N(S)$, $N[S]$ and $d(v)$ rather than
$N_G(v)$, $N_G[v]$, $N_G(S)$, $N_G[S]$ and $d_G(v)$, respectively.

If $G$ is a disjoint union of $k$ copies of a graph $F$, we write $G
= kF$. For a subset $S \subseteq V$, the subgraph induced by $S$ is
denoted by $G[S]$. If $S \subseteq V$, then by $G - S$ we denote the
graph obtained from $G$ by deleting the vertices in the set $S$ (and
all edges incident with vertices in $S$). If $S = \{v\}$, then we
also denote $G - \{v\}$ simply by $G - v$. A \emph{component} in $G$
is a maximal connected subgraph of $G$. If $G$ is a disjoint union of
$k$ copies of a graph $F$, we write $G = kF$. A \emph{star}-\emph{forests} is a forest in which every component is a star.
%For subsets $X,Y \subseteq V$, we denote the set of edges that join a
%vertex of $X$ and a vertex of $Y$ by $[X,Y]$. Thus, $|[X,Y]|$ is the
%number of edges with one end in $X$ and the other end in $Y$.
%%%In particular, $|[X,X]| = m(G[X])$.

\section{Known Results}

Directly from the definition we obtain the following result first observed by Brigham and Dutton~\cite{bd}.

\begin{ob}{\rm (\cite{bd})}
If $G$ is the combined graph of $k \ge 2$ factors, $F_1,F_2,\ldots,F_k$, then
\[
\max_{1 \le i \le k} \gamma(F_i) \le \sd(G) \le \sum_{i=1}^k \gamma(F_i).
\]
 \label{trivial}
\end{ob}

That the lower bound of Observation~\ref{trivial} is sharp, may be seen by taking the $k$ factors, $F_1,F_2,\ldots,F_k$, to be equal. To see that the upper bound of Observation~\ref{trivial} is sharp, let $k \ge 2$ and let $F_1, F_2, \ldots, F_k$ be factors with vertex $V$, where $|V| = n > k$, defined as follows. Let $V = \{v_1,v_2,\ldots,v_n\}$ and let $F_i$ be a star $K_{1,n-1}$ centered at the vertex $v_i$, $1 \le i \le k$. Then, $\{v_1,v_2,\ldots,v_k\}$ is a minimum SD-set of the combined graph $G(F_1,F_2,\ldots,F_k)$, implying that
\[
\sd(F_1,F_2,\ldots,F_k) = \sum_{i = 1}^k  \gamma(F_i) = k.
\]

Brigham and Dutton~\cite{bd} were also the first to observe the following bound.

\begin{ob}{\rm (\cite{bd})}
$\sd(k,\delta,n) \le n - \delta$.
 \label{trivial2}
\end{ob}

The following bounds on $\sd(k,n)$ are established
in~\cite{DaGoHeLa06,DaLa03}.

\begin{thm} The following holds. \\
{\rm (a)} {\rm (\cite{DaLa03})} For $k = 2$, $\sd(k,n) \le 2n/3$, and
this is sharp. \\
{\rm (b)} {\rm (\cite{DaGoHeLa06})} For $k \ge 3$, $\sd(k,n) \le
(2k-3)n/(2k-2)$, and this is sharp for all $k$.
 \label{known_thmA}
\end{thm}

Values of $\sd(k,n)$ in Theorem~\ref{known_thmA} for small $k$ are shown in Table~1.

Caro and Yuster~\cite{CaYu00} considered a combined graph consisting
of $k$ factors $F_1, F_2, \ldots, F_k$. In the language of the
current paper, they were interested in finding a minimum subset $D$
of vertices with the property that the subgraph induced by $D$ is a
connected $r$-dominating set in each of the factors $F_i$, $1 \le i
\le k$, where  $r \le \de = \min \{\, \de(F_i) \mid i=1,2,\ldots,k
\,\}$. As a special consequence of their main result, we have the
following asymptotic result.

\begin{thm}{\rm (\cite{CaYu00})}
Let $F_1, F_2, \ldots, F_k$ be factors on $n$ vertices and let $\de = \min \{\, \de(F_i) \mid i=1,2,\ldots,k \,\}$. If $\delta > 1$ and $\ln \ln \delta > k$, then
\[
\sd(F_1,F_2,\ldots,F_k) \le \left( \frac{(\ln \de)(1 + \oo_\delta(1))}{\de} \right) n.
%n \frac{\ln \delta}{\delta} (1 + \oo_\delta(1)).
\]
 \label{known_thmB}
\end{thm}

Dankelmann and Laskar~\cite{DaLa03} established the following upper bound on the
simultaneous domination number of $k$ factors, depending on the smallest minimum degree of the factors.

\begin{thm} \label{t:general}
Let $F_1, F_2, \ldots, F_k$ be factors on $n$ vertices. Let $\de = \min
\{\, \de(F_i) \mid i=1,2,\ldots,k \,\}$. If $\de \ge 2$ and $k \le e^{\delta + 1}/(\delta + 1)$, then
\[
\sd(F_1,F_2,\ldots,F_k) \le \left( \frac{\ln (\de + 1) + \ln k + 1}{\de + 1} \right) n. \]
\end{thm}

We close this section with a construction showing that the upper bound in Theorem~\ref{known_thmA}(a), which was originally demonstrated by \emph{star}-\emph{forests}, can be realized by trees.
Let $F_1$ and $F_2$ be factors on $n = 3k$ vertices constructed as follows. Let $F_1$ be obtained from the path $u_1u_2 \ldots u_k$ by adding for each $i$, $1 \le i \le k$, two new vertices $v_i$ and $z_i$ and joining $u_i$ to $v_i$ and $z_i$. Further let $F_2$ be obtained from the path $z_1z_2 \ldots z_k$ by adding for each $i$, $1 \le i \le k$, for each $i$, $1 \le i \le k$, add two new vertices $u_i$ and $v_i$ and joining $z_i$ to $u_i$ and $v_i$. We note that both factors $F_1$ and $F_2$ are trees.

Let $D$ be a SD-set of the combined graph $G(F_1,F_2)$. On the one hand, if $u_1 \in D$, then in order to dominate the vertex $v_1$ in $F_2$, we have that at least one of $v_1$ and $z_1$ belong to~$D$. On the other hand, if $u_1 \notin D$, then in order to dominate the vertices $v_1$ and $z_1$ in $F_2$, both $v_1$ and $z_1$ belong to $D$. In both cases, $|D \cap \{u_1,v_1,z_1\}| \ge 2$. Analogously, $|D \cap \{u_i,v_i,z_i\}| \ge 2$ for all $i$, $1 \le i \le k$, implying that $|D| \ge 2k = 2n/3$. Since $D$ was an arbitrary SD-set of $G(F_1,F_2)$, we have that $\sd(F_1,F_2) \ge 2n/3$. Conversely the set $\bigcup_{i=1}^k \{u_i,v_i\}$ is a SD-set of $G(F_1,F_2)$, and so $\sd(F_1,F_2) \le 2n/3$. Consequently, $\sd(F_1,F_2) = 2n/3$ in this case. Further, $\gamma(F_1) = \gamma(F_1) = n/3$. Hence we have the following statement.

\begin{ob}
For $n \equiv 0 \, (\mod \, 3)$, there exist factors $F_1$ and $F_2$ on $n$ vertices, both of which are trees, such that $\sd(F_1,F_2) = 2n/3 = \gamma(F_1) + \gamma(F_2)$.
 \label{t:tree}
\end{ob}

\section{Outline of Paper}

In this paper we continue the study of simultaneous domination in graphs. In Section~\ref{S:gen_bds} we provide general upper bounds on the simultaneous domination number of a combined graph in terms of the generalized vertex cover and independence numbers. Using a hypergraph and probabilistic approach we provide an improvement on the bound of Theorem~\ref{t:general}.
%
%In Section~\ref{S:regular} we focus on combined graphs consisting of regular factors and provide an upper bound on the simultaneous domination number of a combined graph.
%
In Section~\ref{S:clique} we provide general upper bounds on the simultaneous domination number of a combined graph when each factor consists of vertex disjoint union of copies of a clique.
We close in Section~\ref{S:cycle} by studying the case when each factor is a cycle or a disjoint union of cycles.

\section{General Upper Bounds}
\label{S:gen_bds}

%In this section we provide general upper bounds on the %simultaneous domination number of a combined graph.

A vertex and an edge are said to \emph{cover} each other in a graph
$G$ if they are incident in $G$. A \emph{vertex cover} in $G$ is a
set of vertices that covers all the edges of $G$. We remark that a
cover is also called a \emph{transversal} or \emph{hitting set} in
the literature. Thus a vertex cover $T$ has a nonempty intersection
with every edge of $G$. The \emph{vertex covering number} $\TR{G}$ of
$G$ is the minimum cardinality of a vertex cover in $G$. A vertex
cover of size~$\TR{G}$ is called a $\TR{G}$-cover. More generally for
$t \ge 0$ a $t$-\emph{vertex cover} in $G$ is a set of vertices $S$
such that the maximum degree in the graph $G[V \setminus S]$ induced
by the vertices outside $S$ is at most~$t$. The \emph{$t$-vertex
covering number} $\TRt{G}$ of $G$ is the minimum cardinality of a
$t$-vertex cover in $G$. A vertex cover of size~$\TRt{G}$ is called a
$\TRt{G}$-cover.  In particular, we note that a $0$-vertex cover is
simply a vertex cover and that $\TR{G} = \tau_0(G)$.

The independence number $\alpha(G)$ of $G$ is the maximum cardinality
of an independent set of vertices of $G$. More generally, for $k \ge
0$ a $k$-independent set in $G$ is a set of vertices $S$ such that
the maximum degree in the graph $G[S]$ induced by the vertices of $S$
is at most~$k$. The $k$-independence number $\alpha_k(G)$ of $G$ is
the maximum cardinality of a $k$-independent set of vertices of $G$.
In particular, we note that a $0$-independent set is simply an
independent set and that $\alpha(G) = \alpha_0(G)$.

Since the complement of a $t$-vertex cover is a $t$-independent set and conversely, we have the following observation.

\begin{ob}
For a graph $G$ of order~$n$ and an integer $t \ge 0$, we have  $\alpha_t(G) + \TRt{G} = n$.
  \label{cover_indep}
\end{ob}

We recall the following well-known Caro-Wei lower bound on the independence number in terms of the degree sequence of the graph.

\begin{thm}{\rm  (\cite{Ca79,Wei81})} For every graph $G$ of order~$n$,
\[
\alpha(G) \ge \sum_{v \in V(G)} \frac{1}{1 + d_G(v)} \ge \frac{n}{\bard(G) + 1}.
\]
 \label{Caro_Wei}
\end{thm}

We will also need the following recent result by Caro and Hansberg~\cite{CaHa12} who established the following lower bound on the $k$-independence number of a graph.

\begin{thm}{\rm  (\cite{CaHa12})}
For $k \ge 0$ if $G$ is a graph of order~$n$ with average degree $\bard$, then
\[
\alpha_k(G) \ge \left( \frac{ k + 1 }{ \lceil \, \bard \, \rceil + k + 1 } \right) n.
\]
 \label{Caro_Hansberg}
\end{thm}

We begin by establishing the following upper bound on the
simultaneous domination number of a combined graph in terms of the
$t$-vertex cover number and also in terms of the sum of the average
degrees from each factor.

%\newpage
\begin{thm}
Let $F_1, F_2, \ldots, F_k$ be factors on $n$ vertices such that $\delta(F_i) \ge \delta \ge 1$. Let $G  = G(F_1,\ldots,F_k)$ be the combined graph of the factoring $F_1, F_2, \ldots, F_k$, and let $\bard(G) = \bard$ and $\bard(F_i) = \bard_i$ for $i = 1,2,\ldots,k$. Then the following holds. \2 \\
\indent {\rm (a)} $\sd(F_1,F_2,\ldots,F_k) \le \tau_{\delta - 1}(G) = n - \alpha_{\delta - 1}(G)$. \3 \\
\indent {\rm (b)} $\displaystyle{ \sd(F_1,F_2,\ldots,F_k) \le \left( \frac{ \lceil \, \bard \, \rceil }{ \lceil \, \bard \, \rceil + \delta }  \right) \, n }$.
\3 \\
\indent {\rm (c)} If $F_1, F_2, \ldots, F_k$ are regular factors on $n$ vertices each of degree~$\delta$, then
\[
\sd(F_1,F_2,\ldots,F_k) \le \left( \frac{k}{k+1} \right) n.
\]
 \label{thm1}
\end{thm}
\proof Let $G  = G(F_1,\ldots,F_k)$ denote the combined graph of the factoring $F_1, F_2, \ldots, F_k$ and let $G$ have vertex set $V$. By definition of the average degree, we have
\[
\bard = \frac{2m(G)}{n} \le 2 \sum_{i=1}^k \frac{m(F_i)}{n} = \sum_{i=1}^k \frac{2m(F_i)}{n} = \sum_{i=1}^k \bard_i.
\]

(a) Let $S$ be a $\tau_{\delta - 1}(G)$-cover. Hence the graph
$\Delta(G[V \setminus S]) \le \delta - 1$ and $|S| = \tau_{\delta -
1}(G)$. Let $F$ be an arbitrary factor of $G$, and so $F = F_i$ for
some $i \in \{1,2,\ldots,k\}$. Since $\delta(F) \ge \delta$ and since
every vertex in $V \setminus S$ is adjacent to at most $\delta - 1$
other vertices in $V \setminus S$, the set $S$ is a dominating set of
$F$. This is true for each of the $k$ factors in $G(F_1,\ldots,F_k)$.
Therefore, $S$ is a SD-set of $G$, and so $\sd(G) \le |S| =
\tau_{\delta - 1}(G)$. By Observation~\ref{cover_indep}, recall that
$\tau_{\delta - 1}(G) = n - \alpha_{\delta - 1}(G)$.

(b) Since $\delta \ge 1$, we note that $\alpha_{\delta - 1}(G) \ge \alpha_0(G) = \alpha(G)$, implying by Observation~\ref{cover_indep} and Theorem~\ref{Caro_Hansberg} that
\[
\tau_{\delta - 1}(G) = n - \alpha_{\delta - 1}(G) \le n - \left( \frac{ \delta }{ \lceil \, \bard \, \rceil + \delta } \right) n = \left( \frac{ \lceil \, \bard \, \rceil }{ \lceil \, \bard \, \rceil + \delta }  \right) \, n.
\]

The desired result now follows from Part~(a).

(c) Let $F_1, F_2, \ldots, F_k$ be regular factors of degree~$\delta$. Then, $\bard_i = \delta$ for $1 \le i \le k$, and so $\bard \le \sum_{i=1}^k \bard_i = k \delta$. Therefore by Part~(b) above, we have
\[
\sd(F_1,F_2,\ldots,F_k) \le \left( \frac{ \lceil \, \bard \, \rceil }{ \lceil \, \bard \, \rceil + \delta }  \right) \, n \le
\left( \frac{ k \delta }{ (k + 1)\delta  }  \right) \, n = \left( \frac{k}{k+1} \right) n.
\]

This establishes Part~(c), and completes the proof of Theorem~\ref{thm1}.~\qed

\medskip
We next use a hypergraph and probabilistic approach to improve upon a bound already obtained using this approach in~\cite{DaGoHeLa06}.
Let $H$ be a hypergraph. A $k$-\emph{edge} in $H$ is an edge of size~$k$. The \emph{rank} of $H$ is the maximum cardinality among all the edges in $H$.
If all edges have the same cardinality~$k$, the hypergraph is said to be $k$-\emph{uniform}. A subset $T$ of vertices in $H$ is a
\emph{transversal} (also called \emph{vertex cover} or \emph{hitting
set} in many papers) if $T$ has a nonempty intersection with every
edge of $H$. The \emph{transversal number} $\tau(H)$ of $H$ is the
minimum size of a transversal in $H$.
For $r \ge 2$, if $H$ is an $r$-uniform hypergraph with $n$ vertices
and $m$ edges, then it is shown in~\cite{DaGoHeLa06} that $\tau(H)
\le n \le n (\ln (rm/n) + 1)/r$. We improve this bound as follows.

%First we present the following upper bound on the transversal number
%of a uniform hypergraph in terms of the number of vertices $n$ and
%number of edges $m$ and the average degree.

\begin{thm}
For $r \ge 2$, let $H$ be an $r$-uniform hypergraph with $n$ vertices
and $m$ edges and with average degree $d  = rm/n$ and such that
$\delta(H) \ge 1$. Then,
\[
\tau(H) \le \left( 1 - \left( \frac{r-1}{r} \right) \left( \frac{1}{d} \right)^{\frac{1}{r-1}} \right) n \le n( \ln(d) + 1)/r.
\]
\label{hyper_bd_New}
\end{thm}
\proof For $0 \le p \le 1$, choose each vertex in $H$ independently
with probability~$p$. Let $X$ be the set of chosen vertices and let
$Y$ be the set of edges from which no vertex was chosen. Then,
$E(|X|) = np$ and $E(|Y|) = m(1-p^r)$. By linearity of expectation,
we have that $E(|X|+|Y|) = E(|X|) + E(|Y|) = np + m(1-p)^r$. Hence if
we add to $X$ one vertex from each edge in $Y$ we get a transversal
$T$ of $H$ such that $E(|T|) \le np + m(1-p)^r$, implying that
$\tau(H) \le np + m(1-p)^r$. Let $f(p) = np + m(1-p)^r$. This
function is optimized when
\[
p^* = 1 - \left( \frac{1}{d} \right)^{\frac{1}{r-1}},
\]
which is a legitimate value for  $p$ as $d \ge \delta(H) \ge 1$.
Further,
\[
f(p^*) = n - n\left( \frac{1}{d} \right)^{\frac{1}{r-1}} + \left( \frac{nd}{r} \right) \left( \frac{1}{d} \right)^{\frac{r}{r-1}} = \left( 1 - \left( \frac{r-1}{r} \right) \left( \frac{1}{d} \right)^{\frac{1}{r-1}} \right) n.
\]

We also note that $np + m(1-p)^r \le np + me^{-pr}$. Taking $p =
\ln(d)/r = \ln(rm/n)/r  \ge 0$, we get $E(|T|) = E(|X|+|Y|) \le
n\ln(rm/n)/r + n/r  = n( \ln(d) + 1)/r$. Hence the optimal choice of
$p$, namely $p = 1 - \left( \frac{1}{d} \right)^{\frac{1}{r-1}}$,
implies that
\[
\tau(H) \le  \left( 1 - \left( \frac{r-1}{r} \right) \left( \frac{1}{d} \right)^{\frac{1}{r-1}} \right) n \le n( \ln(d) + 1)/r,
\]
which completes the proof of the theorem.~\qed

\medskip
As an application of Theorem~\ref{hyper_bd_New}, we have the
following upper bound on the simultaneous domination number of a combined graph that improves the upper bound of Theorem~\ref{t:general}. For a graph $G$, the \emph{neighborhood hypergraph} of $G$,
denoted by $\NH(G)$, is the hypergraph with vertex set $V(G)$ and
edge set $\{N_G[v] \mid v \in V(G) \}$ consisting of the
closed neighborhoods of vertices in $G$.

\begin{thm}
For $k \ge 2$, if $F_1, F_2, \ldots, F_k$ are
factors on $n$ vertices, each of which has minimum degree at least~$\delta$, then
\[
\sd(F_1,F_2,\ldots,F_k) \le \left( 1 - \left( \frac{\delta}{\delta + 1} \right) \left( \frac{1}{k(\delta + 1)} \right)^{\frac{1}{\delta}} \right) n.
\]
 \label{thm_prob}
\end{thm}
\proof Let $G  = G(F_1,\ldots,F_k)$ denote the combined graph of the
factoring $F_1, F_2, \ldots, F_k$ and let $G$ have vertex set $V$.
Let $\NH(F_i)$ be the neighborhood hypergraph of $F_i$, where $1 \le i \le k$. In particular, we note that $\NH(F_i)$ has vertex set $V$ and rank at least~$\delta + 1$. Let $H_i$ be obtained from $\NH(F_i)$ by shrinking all edges of $\NH(F_i)$, if necessary, to edges of size~$\delta + 1$ (by removing vertices from each edge of size greater than~$\delta + 1$ until the resulting edge size is~$\delta + 1$). Let $H$ be the hypergraph with vertex set $V$ and edge set
$
E(H) = \bigcup_{i=1}^k E(H_i).
$
Then, $H$ is a $(\delta + 1)$-uniform hypergraph with $n(H) = n$ vertices and $m(H) \le kn$ edges. The average degree of $H$ is $d  = (\delta + 1)m(H)/n(H) \le k(\delta + 1)$, implying by Theorem~\ref{hyper_bd_New}, that
\[
\tau(H) \le \left( 1 - \left( \frac{\delta}{\delta + 1} \right) \left( \frac{1}{k(\delta + 1)} \right)^{\frac{1}{\delta}} \right) n.
\]

Every transversal in $H$ is a SD-set in $G$, implying that
$\sd(F_1,F_2,\ldots,F_k) \le \tau(H)$, and the desired result
follows.~\qed

\medskip
Let $f(k,\delta)$ denote the expression on the right-hand side of the inequality in Theorem~\ref{thm_prob}. For small $k$ and small $\delta$, the values of $f(k,\delta)$ are given in Table~3 in the Appendix.

\section{$K_r$-Factors}
\label{S:clique}

As an application of Theorem~\ref{hyper_bd_New}, we have the
following upper bound on the simultaneous domination number of a
combined graph when each factor consists of vertex disjoint union of
copies of $K_r$, for some $r \ge 2$.

\begin{thm}
Let $r$ and $n$ be integers such that $1 \le r \le n$ and $n \equiv 0
\, (\mod \, r)$. For $k \ge 2$, if $F_1, F_2, \ldots, F_k$ are
factors on $n$ vertices, each of which consist of the vertex disjoint
union of $n/r$ copies of $K_r$, then
\[
\sd(F_1,F_2,\ldots,F_k) \le \left( 1 - \left( \frac{r-1}{r} \right) \left( \frac{1}{k} \right)^{\frac{1}{r-1}} \right) n \le n( \ln(k) +1)/r.
\]
 \label{thmB_new}
\end{thm}
\proof Let $G  = G(F_1,\ldots,F_k)$ denote the combined graph of the
factoring $F_1, F_2, \ldots, F_k$ and let $G$ have vertex set $V$.
Let $H$ be the hypergraph with vertex set $V$ and edge set defined as
follows: For every copy of $K_r$ in each of the factors $F_i$, $1 \le
i \le k$, add an $r$-edge in $H$ defined by the vertices of this copy
of $K_r$. The resulting hypergraph $H$ is an $r$-uniform hypergraph
on $n$ vertices with $m \le kn/r$ edges. The average degree of $H$ is
therefore $d  = rm/n \le k$, implying by Theorem~\ref{hyper_bd_New},
that
\[
\tau(H) \le \left( 1 - \left( \frac{r-1}{r} \right) \left( \frac{1}{k} \right)^{\frac{1}{r-1}} \right) n \le n( \ln(k) +1)/r.
\]

Every transversal in $H$ is a SD-set in $G$, implying that
$\sd(F_1,F_2,\ldots,F_k) \le \tau(H)$, and the desired result
follows.~\qed

\medskip
Let $g(k,\delta)$ denote the middle term in the inequality chain in Theorem~\ref{thmB_new}. For small $k$ and small $\delta$, the values of $g(k,\delta)$ are given in Table~4 in the Appendix.

%\medskip
%We note that for $k \ge 4$ and for all $r$ such that $1 \le r \le n$
%and $n \equiv 0 \, (\mod \, r)$, the upper bound of
%Theorem~\ref{thmB_new} is better than that of Theorem~\ref{thmB}.

%For $k = 3$, we remark that for $r \le 6$ the upper bound of
%Theorem~\ref{thmB_new} is smaller than the upper bound of~$2n/r$ from
%Theorem~\ref{thmB} and larger for $r \ge 7$. Hence for $k = 3$ and $r
%\ge 7$ the upper bound of~$2n/r$ in Theorem~\ref{thmB} is better than
%the probabilistic upper bound of Theorem~\ref{thmB_new}, while for $k
%= 3$ and $r \le 6$ the probabilistic bound is better.

Recall that a graph is called \emph{well}-\emph{dominated graph} if
every minimal dominating set in the graph has the same cardinality.
This concept was introduced by Finbow, Hartnell and
Nowakowski~\cite{FiHaNo88}.  We remark that if $v$ is an arbitrary vertex of a well-dominated graph $G$, then the vertex $v$ can be extended to a maximal independent set, which is a minimal dominating set. However, every minimal dominating set in $G$ is a minimum dominating set in $G$ since $G$ is well-dominated. Therefore, every vertex of a well-dominated graph is contained in a minimum dominating set of the graph.

A graph is $1$-\emph{extendable}-\emph{dominated} if every vertex belongs to a minimum dominating set of the graph. We note that every well-dominated graph is a $1$-extendable-dominated graph. However, not every $1$-extendable-dominated graph is well-dominated as may be seen by taking, for example, a cycle $C_6$ or, more generally, a cycle $C_n$, where $n \ge 8$.

\begin{thm}
Let $F$ be a $1$-extendable-dominated graph of order~$r$. Let $n$ be an integer
such that $r \le n$ and $n \equiv 0 \, (\mod \, r)$. If $F_1$ and
$F_2$ are factors on $n$ vertices, each of which consist of the vertex
disjoint union of $n/r$ copies of $F$, then $\sd(F_1,F_2) \le
\frac{1}{r}(2\gamma(F) - 1)n$.
%\[
%\sd(F_1,F_2) \le \left( \frac{2\gamma(F) - 1}{r} \right) n.
%\]
 \label{thmA}
\end{thm}
\proof We construct a bipartite graph $G$ as follows. Let $V_1$ and
$V_2$ be the partite sets of $G$ where for $i \in \{1,2\}$ the
vertices of $V_i$ correspond to the $n/r$ copies of $F$ in $F_i$. An
edge in $G$ joins a vertex $v_1 \in V_1$ and a vertex $v_2 \in V_2$
if and only if the copies of $F$ corresponding to $v_1$ and $v_2$ in
$F_1$ and $F_2$, respectively, have at least one vertex in common. We
observe that $|V_1| = |V_2| = n/r$.

We show that $G$ contains a perfect matching. Let $S$ be a nonempty
subset of vertices of $V_1$. We consider the corresponding $|S|$
vertex disjoint copies of $F$ in $F_1$. These $|S|$ copies of $F$
cover exactly $r|S|$ vertices in $F_1$. But the minimum number of
copies of $F$ in $F_2$ needed to cover these $r|S|$ vertices is at
least $|S|$ since each copy of $F$ covers $r$ vertices. Every vertex
in $V_2$ corresponding to such a copy of $F$ in $F_2$ is joined in
$G$ to at least one vertex of $S$, implying that $|N(S)| \ge |S|$.
Hence by Hall's Matching Theorem, there is a matching in $G$ that
matches $V_1$ to a subset of $V_2$. Since $|V_1| = |V_2|$, such a
matching is a perfect matching in $G$.

Let $M$ be a perfect matching in $G$. For each edge $e \in M$, select
a vertex $v_e$ that is common to the copies of $F$ in $F_1$ and $F_2$
that correspond to the ends of the edge $e$. Since $F$ is a
$1$-extendable-dominated graph, this common vertex $v_e$ extends to minimum
dominating set in both copies of $F$ creating a dominating set of
these two copies with at most $2\gamma(F) - 1$ vertices. Let $D_e$
denote the resulting dominating set of these two copies of $F$. Then
the set $\cup_{e \in M} D_e$ is a SD-set in the combined graph of
$F_1$ and $F_2$, implying that $\sd(F_1,F_2) \le |M| \cdot
(2\gamma(F) - 1) \le (2\gamma(F) - 1)n/r$.~\qed

\medskip
We remark that the bound in Theorem~\ref{thmA} is strictly better than the bound of Theorem~\ref{known_thmA} and Theorem~\ref{thm1}(c) in the case of $k = 2$ when $\gamma(F) < (2r+3)/6$.
As a consequence of Theorem~\ref{thmA}, we have the following results.

\begin{thm}
Let $r$ and $n$ be integers such that $1 \le r \le n$ and $n \equiv 0
\, (\mod \, r)$. If $F_1$ and $F_2$ are factors on $n$ vertices, each of
which consist of the vertex disjoint union of $n/r$ copies of $K_r$,
then $\sd(F_1,F_2) = n/r$.
 \label{t:twoF}
\end{thm}
\proof We note that $K_r$ is a well-dominated graph. Further,
$\gamma(K_r) = 1$. Applying Theorem~\ref{thmA} with the graph $F =
K_r$, we have that $\sd(F_1,F_2) \le n/r$. By
Observation~\ref{trivial}(a), we know that $\sd(F_1,F_2) \ge
\gamma(F_1) = n/r$. Consequently, $\sd(F_1,F_2) = n/r$.~\qed

\begin{cor}
Let $r$ and $n$ be integers such that $1 \le r \le n$ and $n \equiv 0
\, (\mod \, r)$. If $F_1$ and $F_2$ are factors on $n$ vertices, each of
which contain a spanning subgraph that is the vertex disjoint union
of $n/r$ copies of $K_r$, then $\sd(F_1,F_2) \le n/r$.
 \label{cor2}
\end{cor}

As an immediate consequence of Corollary~\ref{cor2} and Observation~\ref{trivial}, we have the following observation.

\begin{cor}
For $n$ even, if $F_1$ and $F_2$ are factors on $n$ vertices both having a $1$-factor, then $\sd(F_1,F_2) \le n/2$. Further, if $\max \{\gamma(F_1), \gamma(F_2)\} = n/2$, then $\sd(F_1,F_2) = n/2$.
 \label{c:1factor}
\end{cor}

We next extend the result of Theorem~\ref{t:twoF} to more than two factors.

\begin{thm}
Let $r$ and $n$ be integers such that $1 \le r \le n$ and $n \equiv 0
\, (\mod \, r)$. For $k \ge 2$, if $F_1, F_2, \ldots, F_k$ are
factors on $n$ vertices, each of which consist of the vertex disjoint
union of $n/r$ copies of $K_r$, then
\[
\sd(F_1,F_2,\ldots,F_k)
\le \left( 1 - \left( \frac{r-1}{r} \right)^{k-1} \right) n.
\]
 \label{thm_kfactors}
\end{thm}
\proof We proceed by induction on $k \ge 2$. The base case when $k =
2$ follows from Theorem~\ref{t:twoF}. Assume, then, that $k \ge 3$
and that the result holds for $k'$ factors, each of which consist of
the vertex disjoint union of $n/r$ copies of $K_r$, where $2 \le k' <
k$. Let $F_1, F_2, \ldots, F_k$ be factors on $n$ vertices, each of
which consist of the vertex disjoint union of $n/r$ copies of $K_r$.
First we consider the combined graph $G(F_1,F_2,\ldots,F_{k-1})$ with
only $F_1, F_2, \ldots, F_{k-1}$ as factors. Let $D$ be a
$\sd(F_1,F_2,\ldots,F_{k-1})$-set in $G(F_1,F_2,\ldots,F_{k-1})$, and
so $|D| = \sd(F_1,F_2,\ldots,F_{k-1})$. By the inductive hypothesis,

\[
|D| \le \left( 1 - \left( \frac{r-1}{r} \right)^{k-2} \right) n.
\]

We now consider the combined graph $G(F_1,F_2,\ldots,F_k)$. Since
each copy of $K_r$ in $F_k$ can have at most $r$ vertices from $D$,
the set $D$ must dominate at least $|D|/r$ copies of $K_r$ from
$F_k$. Therefore in $F_k$ there remains at most $n/r - |D|/r$ copies
of $F_k$ that are not dominated by $D$. We now extend the set $D$ to
an SD-set of $G(F_1,F_2,\ldots,F_k)$ by adding to it one vertex from
each non-dominated copy of $K_r$ of $F_k$. Hence,

\[
\begin{array}{lcl} \2
\sd(F_1,F_2,\ldots,F_k) & \le & \displaystyle{ |D| + \frac{n-|D|}{r}  } \\ \2
& = & \displaystyle{ \frac{1}{r} ( n + (r-1)|D|) } \\ \3
& \le & \displaystyle{ \frac{1}{r} \left( n + (r-1)
\left( 1 - \left( \frac{r-1}{r} \right)^{k-2} \right) n  \right) } \\ \3
& \le & \displaystyle{ \frac{1}{r} \left( r - (r-1)
\left( \frac{r-1}{r} \right)^{k-2} \right) n   } \\ \3
& = & \displaystyle{ \left( 1 - \left( \frac{r-1}{r} \right)^{k-1} \right) n },
\end{array}
\]

\noindent completing the proof of the theorem.~\qed

\medskip
We remark that the bound in Theorem~\ref{thm_kfactors} is strictly better than the bounds of  Theorem~\ref{known_thmA}, Theorem~\ref{thm1}(c) and Theorem~\ref{thmB_new} when $k = 3$ and for all $r \ge 3$. In particular, we remark that when $k = 3$ and $r \ge 3$, the bound in Theorem~\ref{thm_kfactors} is strictly better than the bound of Theorem~\ref{thmB_new} if
\[
1 - \left( \frac{r-1}{r} \right)^{2}  <
1 - \left( \frac{r-1}{r} \right) \left( \frac{1}{3} \right)^{\frac{1}{r-1}},
\]
or, equivalently, if
\[
\frac{1}{3} < \left( \frac{r-1}{r} \right)^{r-1}.
\]

Since $\left( \frac{r-1}{r} \right)^{r-1}$ attains the value~$4/9$ when $r = 3$ and is a decreasing function in $r$ approaching $0.367879$ as $r \rightarrow \infty$, the above inequality holds.
In the special case in Theorem~\ref{thm_kfactors} when $k = 3$, we have the following result.

\begin{cor}
Let $r$ and $n$ be integers such that $1 \le r \le n$ and $n \equiv 0 \, (\mod \, r)$.  If $F_1, F_2, F_3$ are factors on $n$ vertices, each of which consist of the vertex disjoint union of $n/r$ copies of $K_r$,
then
\[
\sd(F_1,F_2,F_3) \le \left( \frac{2r-1}{r^2} \right) \, n.
\]
 \label{c:k3r3}
\end{cor}

Using Corollary~\ref{c:k3r3}, the upper bound of
Theorem~\ref{thm_kfactors} can be improved slightly as follows.

\begin{thm}
Let $r$ and $n$ be integers such that $1 \le r \le n$ and $n \equiv 0
\, (\mod \, r)$. For $k \ge 2$, if $F_1, F_2, \ldots, F_k$ are
factors on $n$ vertices, each of which consist of the vertex disjoint
union of $n/r$ copies of $K_r$, then

\[
\sd(F_1,F_2,\ldots,F_k) \le \left\{
\begin{array}{cl} \2
 \displaystyle{ \left( \frac{k}{2r} \right) n } & \mbox{if $k$ is even}   \\
 \displaystyle{ \left( \frac{r(k+1) - 2}{2r^2} \right) n } & \mbox{if $k$ is odd}.
\end{array}
\right.
\]
%$\sd(F_1,F_2,\ldots,F_k) \le kn/2r$.
 \label{thm_kfactors_B}
\end{thm}
\proof Suppose first that $k$ is even. Consider the combined graph
$G(F_{2i-1},F_{2i})$ with only $F_{2i-1}$ and $F_{2i}$ as factors,
where $1 \le i \le k/2$. For each such~$i$, let $D_i$ be a
$\sd(F_{2i-1},F_{2i})$-set in $G(F_{2i-1},F_{2i})$ and note that by
Theorem~\ref{t:twoF}, we have $|D_i| = n/r$. Let $D =
\bigcup_{i=1}^{k/2} D_i$. Then the set $D$ is a SD-set of
$G(F_1,F_2,\ldots,F_k)$, implying that $\sd(F_1,F_2,\ldots,F_k) \le
|D| \le kn/2r$.

Suppose next that $k$ is odd. Let $D_1$ be a $\sd(F_1,F_2,F_3)$-set
in the combined graph $G(F_1,F_2,F_3)$ with only $F_1,F_2,F_3$ as
factors. By Corollary~\ref{c:k3r3}, we have $|D_1| \le (2r-1)n/r^2$.
For $i$ with $2 \le i \le (k-1)/2$, consider the combined graph
$G(F_{2i},F_{2i+1})$ with only $F_{2i}$ and $F_{2i+1}$ as factors and
let $D_i$ be a $\sd(F_{2i},F_{2i+1})$-set in $G(F_{2i},F_{2i+1})$. By
Theorem~\ref{t:twoF}, we have $|D_i| = n/r$ for $2 \le i \le
(k-1)/2$. Let $D = \bigcup_{i=1}^{(k-1)/2} D_i$. Then the set $D$ is
a SD-set of $G(F_1,F_2,\ldots,F_k)$, implying that
\[
\sd(F_1,F_2,\ldots,F_k) \le |D| \le
\left( \frac{2r-1}{r^2} \right) \, n + \left( \frac{k-3}{2r} \right) \, n
= \left( \frac{r(k+1) - 2}{2r^2} \right) n,
\]
which established the desired upper bound in this case when $k$ is
odd.~\qed

\medskip
We remark that the bound in Theorem~\ref{thm_kfactors_B} is strictly better than the bounds of  Theorem~\ref{known_thmA} and Theorem~\ref{thm1}(c) for $r \ge 3$. Further the bound in Theorem~\ref{thm_kfactors_B} is strictly better than the bound of Theorem~\ref{thmB_new} for $r \ge 4$.

We close this section by
%
%\subsection{$K_2$-Factors}
%\label{S:K2}
%
%\medskip
considering the special case when every factor in the
combined graph is the disjoint union of copies of $K_2$. If $G$ is a
graph of even order and if $F$ is a $1$-regular spanning subgraph of
$G$, we call $F$ a $1$-\emph{factor} of $G$. Hence if $F$ is a
$1$-factor of a graph $G$ of order~$n$, then $F = \frac{n}{2}K_2$ and
the edges of $F$ form a perfect matching in $G$.

\begin{thm}
%Let $n$ be an even integer.
For $k \ge 2$ and $n$ even, if $F_1, F_2, \ldots, F_k$ are
$1$-factors on $n$ vertices, then

\[
\sd(F_1,F_2,\ldots,F_k) \le \left\{
\begin{array}{cl} \2
 \displaystyle{ \left( \frac{k-1}{k} \right) n } & \mbox{if $k$ is even}   \\
 \displaystyle{ \left( \frac{k}{k+1} \right) n } & \mbox{if $k$ is odd}.
\end{array}
\right.
\]
and these bounds are sharp.
 \label{thmC}
\end{thm}
\proof Let $G  = G(F_1,\ldots,F_k)$ denote the combined graph of the
factoring $F_1, F_2, \ldots, F_k$ and let $G$ have vertex set $V$.
Then, $\Delta(G) \le k$. By Brook's Coloring Theorem, $\chi(G) \le
k+1$ with equality if and only if $G$ has a component isomorphic to
$K_{k+1}$ or a component that is an odd cycle and $k = 2$. %

We show that every
component of $G$ has even order. Suppose to the contrary that there
is a component, $F$, in $G$ of odd order. For each vertex $v$ in
$V(F)$, let $v'$ be its neighbor in $F_1 = \frac{n}{2}K_2$ and let $S
= \cup_{v \in V(F)} \{v,v'\}$. Then, $V(F) = S$. However, $|S|$ is
even, while $|V(F)|$ is odd, a contradiction. Therefore, every
component of $G$ has even order. In particular, no component of $G$
is an odd cycle.

If $k$ is odd, then by Theorem~\ref{thm1}(c),
$\sd(F_1,F_2,\ldots,F_k) \le kn/(k+1)$, as desired. If $k$ is even,
then no component of $G$ is isomorphic to $K_{k+1}$, implying that
$\chi(G) \le k$. This in turn implies that $\alpha(G) \ge n/\chi(G) =
n/k$, and so, by Observation~\ref{cover_indep} and
Theorem~\ref{thm1}(a) we have that $\sd(F_1,F_2,\ldots,F_k) \le
\tau(G) = n - \alpha(G) \le (k-1)n/k$, as desired.

That these bounds are sharp may be seen as follows. For $k$ odd, take
$n \equiv 0 \, (\mod \, k+1)$. Then the $1$-factors $F_1, F_2,
\ldots, F_k$ of $K_n$ can be chosen so that the combined graph $G$
consists of the disjoint union of $n/(k+1)$ copies of $K_{k+1}$. Let
$S$ be an SD-set in $G$ of minimum cardinality and let $F$ be an
arbitrary copy of $K_{k+1}$ in $G$. If $|S \cap V(F)| \le k-1$, then
there would be two vertices, $u$ and $v$, in $F$ that do not belong
to $S$. However the edge $uv$ belongs to one of the factor of $G$,
implying that in such a $1$-factor neither $u$ nor $v$ is dominated
by $S$, a contradiction. Hence, $|S \cap V(F)| \ge k$. This is true
for every copy of $K_{k+1}$ in $G$. Therefore,
$\sd(F_1,F_2,\ldots,F_k) = |S| \ge kn/(k+1)$. As shown earlier,
$\sd(F_1,F_2,\ldots,F_k) \le kn/(k+1)$. Consequently,
$\sd(F_1,F_2,\ldots,F_k) = kn/(k+1)$.

For $k$ even, we simply take $F_{k-1} = F_k$, and note that in this
case $\sd(F_1,F_2,\ldots,F_k) = \sd(F_1,F_2,\ldots,F_{k-1})$. Since
$k-1$ is odd, the construction in the previous paragraph shows that
the $1$-factors $F_1, F_2, \ldots, F_{k-1}$ of $K_n$ can be chosen so
that the combined graph $G$ satisfies $\sd(F_1,F_2,\ldots,F_k) =
(k-1)n/k$.~\qed

\medskip
We remark that the bound in Theorem~\ref{thmC} is better than the bound of Theorem~\ref{thmB_new} always, better than the bound of Theorem~\ref{known_thmA} for $k \ge 2$, and better than the bound of Theorem~\ref{thm1}(c) for $k$ even.

\section{Cycle Factors}
\label{S:cycle}

In this section, we consider the case when each factor is a cycle or a disjoint union of cycles. As a consequence of Corollary~\ref{cor2}, we have the following upper
bound on the simultaneous domination number of a combined graph with
two factors, both of which are cycles or paths.

\begin{thm} The following holds. \\
\indent {\rm (a)} For $n \equiv 0 \, (\mod \, 2)$ and $n \ge 4$,
$\sd(C_n,C_n) \le n/2$ and $\sd(P_n,P_n) \le n/2$. \\
\indent {\rm (b)} For $n \equiv 1 \, (\mod \, 2)$ and $n \ge 5$, $\sd(C_n,C_n) \le
(n+1)/2$.
 \label{t:cycleF}
\end{thm}
\proof (a) For $n \equiv 0 \, (\mod \, 2)$ and $n \ge 4$, both the
cycle $C_n$ and the path $P_n$ contains a spanning subgraph that is
the vertex disjoint union of $n/2$ copies of $K_2$, and so by
Corollary~\ref{cor2}, we have that $\sd(C_n,C_n) \le n/2$ and
$\sd(P_n,P_n) \le n/2$.

(b) For $n \equiv 1 \, (\mod \, 2)$ and $n \ge 3$, let $v$ be an
arbitrary vertex in the cycle $C_n$. Deleting the vertex $v$ from the
cycle, we produce a path $P_{n-1}$, where $n - 1 \equiv 0 \, (\mod \,
2)$. Applying Part~(a), we have that $\sd(P_{n-1},P_{n-1}) \le
(n-1)/2$. Adding the deleted vertex $v$ to a minimum SD-set in the
combined graph with the two paths $P_{n-1}$ as factors, we produce a
SD-set in the original combined graph with the two cycles $C_n$ as
factors of cardinality~$\sd(P_{n-1},P_{n-1}) + 1 \le (n+1)/2$.~\qed

\medskip
For generally, we can establish the following upper bound on the
simultaneous domination number of a combined graph with $k \ge 2$
factors, each of which is a cycle. For simplicity, we restrict the
number of vertices to be congruent to zero modulo~$6$.

\begin{thm}
For $k \ge 2$ and $n \equiv 0 \, (\mod \, 6)$, let $F_1, F_2, \ldots,
F_k$ be factors on $n$ vertices, each of which is isomorphic to a
cycle $C_n$. Then,
\[
\sd(F_1,F_2,\ldots,F_k) \le \left(  1 - \frac{1}{2} \left( \frac{2}{3} \right)^{k-2}  \right) \, n.
\]
 \label{t:cycle_Gen}
\end{thm}
\proof We proceed by induction on $k \ge 2$. The base case when $k =
2$ follows from Theorem~\ref{t:cycleF}(a). Assume, then, that $k \ge
3$ and that the result holds for $k'$ factors, each of which is
isomorphic to a cycle $C_n$, where $2 \le k' < k$. Let $F_1, F_2,
\ldots, F_k$ be factors on $n$ vertices, each of which is isomorphic
to a cycle $C_n$.
First we consider the combined graph $G(F_1,F_2,\ldots,F_{k-1})$ with
only $F_1, F_2, \ldots, F_{k-1}$ as factors. Let $D$ be a
$\sd(F_1,F_2,\ldots,F_{k-1})$-set in $G(F_1,F_2,\ldots,F_{k-1})$, and
so $|D| = \sd(F_1,F_2,\ldots,F_{k-1})$. By the inductive hypothesis,

\[
|D| \le \left(  1 - \frac{1}{2} \left( \frac{2}{3} \right)^{k-3}  \right) \, n.
\]

We now consider the combined graph $G(F_1,F_2,\ldots,F_k)$. Let $F_k$
be the cycle $v_1v_2 \ldots v_n v_1$. For $i = 1,2,3$, let $D_i =
\{v_j \mid j \equiv i \, (\mod \, 3) \}$. We note that for $i \in
\{1,2,3\}$, each set $D_i$ is a dominating set in $F_k$ and $|D_i| =
n/3$. We now extend the set $D$ to a SD-set of
$G(F_1,F_2,\ldots,F_k)$ as follows. Renaming vertices, if necessary,
we may assume that
\[
|D \cap D_1| = \max_{1 \le i \le 3} |D \cap D_i|.
\]

Thus, $|D| = \sum_{i=1}^3 |D \cap D_i| \le 3|D \cap D_1|$, or,
equivalently, $|D \cap D_1| \ge |D|/3$. Let $S $ be the set of
vertices in $D_1$ that do belong to $D$. Then, $S = D_1 \setminus D$
and $|S| = |D_1| - |D \cap D_1| \le n/3 - |D|/3$. Since $D_1
\subseteq D \cup S$ and $D_1$ is a dominating set of $F_k$, the set
$D \cup S$ is a dominating set of $F_k$. Since $D$ is a DS-set of
$G(F_1,F_2,\ldots,F_{k-1})$, the set $D$ is a dominating set in $F_i$
for $1 \le i \le k-1$. Hence, $D \cup S$ is a SD-set of
$G(F_1,F_2,\ldots,F_k)$, implying that

\[
\begin{array}{lcl} \2
\sd(F_1,F_2,\ldots,F_k) & \le & |D| + |S| \\ \2
& \le & |D| + \frac{n - |D|}{3} \\ \2
& \le & \frac{n + 2|D|}{3} \\ \2
& \le & \frac{1}{3} \left( n + 2 \left(  1 - \frac{1}{2} \left( \frac{2}{3} \right)^{k-3}  \right)  \, n  \right) \\ \2
& = &  \left(  1 - \frac{1}{2} \left( \frac{2}{3} \right)^{k-2}  \right) \, n,
\end{array}
\]

\noindent completing the proof of the theorem.~\qed

\medskip
We remark that Theorem~\ref{t:cycle_Gen} is better than Theorem~\ref{thmC} when $k = 3$, since in this case the upper bound of Theorem~\ref{t:cycle_Gen} is~$2n/3$ while that of Theorem~\ref{thmC} is~$3n/4$.

\subsection{$C_4$-Factors}

We consider here the case when every factor in the combined graph is
the disjoint union of copies of a $4$-cycle. As a consequence of Corollary~\ref{c:1factor}, we have the following result.

\begin{thm}
For $n \equiv 0 \, (\mod \, 4)$, let $F_1$ and $F_2$ be factors on
$n$ vertices, both of which are isomorphic to $\frac{n}{4}C_4$. Then,
$\sd(F_1,F_2) = n/2$.
 \label{t:cycle_2C4}
\end{thm}
\proof We observe that $F_1$ and $F_2$ are factors on $n$ vertices both having a $1$-factor. Further, each of the $n/4$ copies of $C_4$ in $F_1$ need two vertices to dominate that copy of $C_4$, implying that $\gamma(F_1) \ge n/2$. The desired result now follows from Corollary~\ref{c:1factor}.~\qed

\begin{thm}
For $n \equiv 0 \, (\mod \, 4)$, let $F_1, F_2, F_3$ be factors on
$n$ vertices, each of which is isomorphic to $\frac{n}{4}C_4$. Then,
$\sd(F_1,F_2,F_3) \le 3n/4$.
 \label{t:cycle_3C4}
\end{thm}
\proof First we consider the combined graph $G(F_1,F_2)$ with only
$F_1$ and $F_2$ as factors. Let $D$ be a $\sd(F_1,F_2)$-set in
$G(F_1,F_2)$. By Theorem~\ref{t:cycle_2C4}, $|D| = n/2$. We next
consider the factor $F_3$. For $0 \le i \le 4$, let $n_i$ denote the
number of copies of $C_4$ in $F_3$ that contain exactly~$i$ vertices
in the set $D$. Counting the number of vertices not in $D$, we have
that
\[
\frac{n}{2} = n - |D| = \sum_{i=0}^4 (4-i)n_i \ge 4n_0 + 3n_1,
\]
implying that $2n_0 + n_1 \le 2n_0 + 3n_1/2 \le n/4$. We now extend
the set $D$ to a SD-set of $G(F_1,F_2,F_3)$ as follows. From each
copy of $C_4$ in $F_3$ that contains exactly one vertex in $D$, we
add to $D$ the vertex that is not adjacent in $F_3$ to a vertex of
$D$. From each copy of $C_4$ in $F_3$ that contains no vertex in $D$,
we add any two vertices to $D$. The resulting set is a SD-set of
$G(F_1,F_2,F_3)$, implying that $\sd(F_1,F_2,F_3) \le |D| + 2n_0 +
n_1 \le n/2 + n/4 = 3n/4$.~\qed

\medskip
We remark that the bound in Theorem~\ref{t:cycle_2C4} is strictly better than the bounds of  Theorem~\ref{known_thmA} and Theorem~\ref{thm1}(c) when $k = 2$. The bound in Theorem~\ref{t:cycle_3C4}, namely $3n/4$, is better than the general probabilistic bound of Theorem~\ref{thm_prob}, namely $f(3,2)n = 7n/9$ (see Table~3).

\subsection{$C_5$-Factors}

We consider here the case when every factor in the combined graph is the disjoint union of copies of a $5$-cycle.

\begin{thm}
For $n \equiv 0 \, (\mod \, 5)$ and $k \ge 2$, let $F_1, F_2, \ldots, F_k$ be factors on $n$ vertices, each of which is isomorphic to $\frac{n}{5}C_5$. Then, $\sd(F_1,F_2) \le 3n/5$ and this bound is sharp. Further, for $k \ge 3$,
\[
\sd(F_1,F_2,\ldots,F_k) \le \left( \frac{3}{5} +  \frac{2}{5} \left( 1 - \left( \frac{3}{5} \right)^{k-2} \right) \right) \, n.
\]
 \label{t:cycle_rC5}
\end{thm}
\proof  We proceed by induction on $k \ge 2$. Let $F_1$ and $F_2$ be factors on $n$ vertices, where both $F_1$ and $F_2$ consist of the vertex-disjoint union of $n/5$ copies of $C_5$. Since the $5$-cycle $C_5$ is well-dominated, we have by Theorem~\ref{thmA} that $\sd(F_1,F_2) \le \frac{1}{5}(2\gamma(C_5) - 1)n = 3n/5$. This establishes the base case when $k = 2$. Assume, then, that $k \ge 3$ and that the result holds for $k'$ factors, each of which consist of the vertex disjoint union of $n/5$ copies of $C_5$, where $2 \le k' < k$. Let $F_1, F_2, \ldots, F_k$ be factors on $n$ vertices, each of which is isomorphic to $\frac{n}{5}C_5$.
First we consider the combined graph $G(F_1,F_2,\ldots,F_{k-1})$ with only $F_1, F_2, \ldots, F_{k-1}$ as factors. Let $D'$ be a $\sd(F_1,F_2,\ldots,F_{k-1})$-set in $G(F_1,F_2,\ldots,F_{k-1})$, and so $|D'| = \sd(F_1,F_2,\ldots,F_{k-1})$.
By the inductive hypothesis, $|D'| \le 3n/5$ if $k = 3$, while for $k \ge 4$, we have

\[
|D'| \le \left( \frac{3}{5} +  \frac{2}{5} \left( 1 - \left( \frac{3}{5} \right)^{k-3} \right) \right) \, n.
\]

We add vertices to $D'$, if necessary, until the cardinality of the resulting superset $D$ is either~$3n/5$ if $k = 3$ or is precisely the expression on the right-hand side of the above inequality if $k \ge 4$. Since $D'$ is a SD-set of $G(F_1,F_2,\ldots,F_{k-1})$, so too is the set $D$. We now consider the combined graph $G(F_1,F_2,\ldots,F_k)$. For $0 \le i \le 5$, let $n_i$ denote the number of copies of $C_5$ in $F_k$ that contain exactly~$i$ vertices in the set $D$. Counting the number of vertices not in $D$, we have that

\[
\frac{2}{5} \left( \frac{3}{5} \right)^{k-3} n = n - |D| = \sum_{i=0}^5 (5-i)n_i \ge 5n_0 + 4n_1 + 3n_2 \ge 5n_0 + 5( n_1 + n_2 )/2,
\]
implying that
\[
2n_0 + n_1 + n_2 \le \frac{4}{25} \left( \frac{3}{5} \right)^{k-3} n.
\]

We now extend the set $D$ to a SD-set of $G(F_1,F_2,\ldots,F_k)$ as follows. From each copy of $C_5$ in $F_k$ that contains no vertex of $D$, we add two vertices that dominate that copy of $C_5$. From each copy of $C_5$ in $F_k$ that contains one or two vertices of $D$, we select one such vertex of $D$ and we add to $D$ a vertex from that copy of $C_5$ that is not adjacent in $F_k$ to that selected vertex. The resulting set is a SD-set of
$G(F_1,F_2,\ldots,F_k)$, implying that
\[
\sd(F_1,F_2,\ldots,F_k) \le |D| + 2n_0 + n_1 + n_2.
\]
If $k = 3$, then
\[
\sd(F_1,F_2,\ldots,F_k) \le \frac{3n}{5} + \frac{4n}{25} = \left( \frac{3}{5} +  \frac{2}{5} \left( 1 - \left( \frac{3}{5} \right)^{k-2}
 \right) \right) \, n.
\]
If $k \ge 4$, then

\[
\begin{array}{lcl} \3
\sd(F_1,F_2,\ldots,F_k) & \le & \displaystyle{ \left( \frac{3}{5} +  \frac{2}{5} \left( 1 - \left( \frac{3}{5} \right)^{k-3} \right) \right) \, n + \frac{4}{25} \left( \frac{3}{5} \right)^{k-3} n} \\ \3
& = & \displaystyle{ \left( \frac{3}{5} +  \frac{2}{5} \left( 1 - \left( \frac{3}{5} \right)^{k-3} + \frac{2}{5}
\left( \frac{3}{5} \right)^{k-3} \right) \right) \, n } \\ \3
& = & \displaystyle{ \left( \frac{3}{5} +  \frac{2}{5} \left( 1 - \frac{3}{5} \left( \frac{3}{5} \right)^{k-3} \right) \right) \, n } \\ \3
& = & \displaystyle{
\left( \frac{3}{5} +  \frac{2}{5} \left( 1 - \left( \frac{3}{5} \right)^{k-2} \right) \right) \, n. }
\end{array}
\]

\noindent completing the proof of the upper bound of the theorem. That the bound is sharp when $k \ge 2$, may be seen as follows. For $r \ge 1$, let $G = rK_5$ be the disjoint union of $r$ copies of $K_5$ and let $G$ have order~$n$. Then there exists two edge-disjoint spanning subgraphs, $F_1$ and $F_2$, of $G$ both of which are isomorphic to the disjoint union of $r$ copies of $C_5$. In order to simultaneously dominate the copies of $C_5$ in $F_1$ and $F_2$ corresponding to a copy of $K_5$ in $G$ at least three vertices are needed, implying that $\sd(F_1,F_2) \ge 3r = 3n/5$. By Theorem~\ref{t:cycle_rC5}, $\sd(F_1,F_2) \le 3n/5$. Consequently, $\sd(F_1,F_2) = 3n/5$ in this case.~\qed

\medskip
We remark that the bound in Theorem~\ref{t:cycle_rC5} is strictly better than the bounds of Theorem~\ref{known_thmA} and Theorem~\ref{thm1}(c) when $k = 2$.
Theorem~\ref{t:cycle_rC5} (when $k = 2$) implies the following result.

\begin{thm}
$\sd(2,2,n) \ge 3n/5$.
 \label{t:22n}
\end{thm}

\section{Open Questions and Conjectures}

Recall that in Theorem~\ref{t:22n}, we established that $\sd(2,2,n)
\ge 3n/5$. The following conjecture was posed by Dankelmann and
Laskar~\cite{DaLa03}, albeit using different notation.

\begin{conj}
$\sd(2,2,n) = 3n/5$.
 \label{conj1}
\end{conj}

By Theorem~\ref{t:cycle_rC5}, if Conjecture~\ref{conj1} is true, then it suffices to prove the following statement: If $F_1$ and $F_2$
are factors on $n$ vertices both having minimum degree at least~$2$,
then $\sd(F_1,F_2) \le 3n/5$.

Recall that in Theorem~\ref{t:cycleF}, for $n \equiv 0 \, (\mod \,
2)$ and $n \ge 4$, we show that $\sd(C_n,C_n) \le n/2$ and
$\sd(P_n,P_n) \le n/2$. Further for $n \equiv 1 \, (\mod \, 2)$ and
$n \ge 5$, $\sd(C_n,C_n) \le (n+1)/2$. We pose the following problem.

\begin{prob}
For all $n \ge 4$, determine the exact value of $\sd(C_n,C_n)$ and
$\sd(P_n,P_n)$.
\end{prob}

Recall by Corollary~\ref{c:1factor} that if $F_1$ and $F_2$ are factors on $n$ vertices both having a $1$-factor, then $\sd(F_1,F_2) \le n/2$. Further, if $\max \{\gamma(F_1), \gamma(F_2)\} = n/2$, then $\sd(F_1,F_2) = n/2$. We close with the following problem that we have yet to settle.

\begin{prob}
Characterize the connected factors $F_1$ and $F_2$ on $n$ vertices that have a $1$-factor and satisfy $\sd(F_1,F_2) = n/2$.
 \label{prob3}
\end{prob}

%For subsets $X,Y \subseteq V(G)$ of a graph $G$, we denote %the set of edges that join a vertex of $X$ and a vertex of %$Y$ by $[X,Y]$.
For $n$ even, let $\cG$ be the family of graphs $G$ whose vertex set can be partitioned into two sets $X$ and $Y$ such that $|X| = |Y| = n/2$, the set $[X,Y]$ of edges that join a vertex of $X$ and a vertex of $Y$ is a $1$-factor in $G$, the set $X$ is independent, and the subgraph $G[Y]$ is connected. By construction, every graph in the family $\cG$ is connected, has a $1$-factor and has domination number one-half its order. Therefore by Corollary~\ref{c:1factor}, we observe that if $F_1$ and $F_2$ are factors on $n$ vertices that belong to the family~$\cG$, then $\sd(F_1,F_2) = n/2$. However we have yet to provide a characterization of all factors $F_1$ and $F_2$ that meet the requirements of Problem~\ref{prob3}.

\newpage
%\medskip

\newpage

\begin{center}
\underline{\textbf{APPENDIX:}}
\end{center}

\hspace*{1cm} \underline{\hspace{12cm}}

\[
\begin{array}{c|cccccc}
k        &    2  &  3  &  4 & 5 & 6 & 7 \\ \hline
& & & & & & \\
%\sd(k,n)  &   \frac{2}{3}   &  \frac{3}{4}  &  \frac{5}{6} %& \frac{7}{8} & \frac{9}{10} & \frac{11}{12}  \\
\sd(k,n)  &   \displaystyle{ \frac{2}{3} }  &  \displaystyle{ \frac{3}{4} } &  \displaystyle{ \frac{5}{6} } & \displaystyle{ \frac{7}{8} } & \displaystyle{ \frac{9}{10} } & \displaystyle{ \frac{11}{12} } \\
\end{array}
\]
\begin{center}
\textbf{Table~1.} Upper bounds on $\sd(k,n)$ in Theorem~\ref{known_thmA} for small $k$.
\end{center}

\hspace*{1cm} \underline{\hspace{12cm}}

\[
\begin{array}{c|cccccc}
k        &    2  &  3  &  4 & 5 & 6 & 7 \\ \hline
& & & & & & \\
%\sd(k,n)  &   \frac{2}{3}   &  \frac{3}{4}  &  \frac{5}{6} %& \frac{7}{8} & \frac{9}{10} & \frac{11}{12}  \\
\sd(k,n)  &   \displaystyle{ \frac{2}{3} }  &  \displaystyle{ \frac{3}{4} } &  \displaystyle{ \frac{4}{5} } & \displaystyle{ \frac{5}{6} } & \displaystyle{ \frac{6}{7} } & \displaystyle{ \frac{7}{8} } \\
\end{array}
\]
\begin{center}
\textbf{Table~2.} Upper bounds on $\sd(k,n)$ in Theorem~\ref{thm1}(c) for small $k$.
\end{center}

\hspace*{1cm} \underline{\hspace{12cm}}

\[
\begin{array}{cc|cccc}
& & & $\hspace*{1cm} k$ & & \\
    &     &    $2$   &    $3$   &     $4$ & $5$ \\ \hline
    & $1$ & $0.8750$ & $0.9167$ & $0.9375$ & $0.9500$ \\
    & $2$ & $0.7278$ & $0.7777$ & $0.8075$ & $0.8278$ \\
$r$ & $3$ & $0.6250$ & $0.6724$ & $0.7023$ & $0.7237$ \\
    & $4$ & $0.5501$ & $0.5935$ &  $0.6217$ & $0.6432$ \\
    & $5$ & $0.4930$ & $0.5325$ &  $0.5586$ & $0.5779$ \\
\end{array}
\]
\begin{center}
\textbf{Table~3.} Approximate values of $f(k,\delta)$ in Theorem~\ref{thm_prob} for small $k$ and $\delta$.
\end{center}

\hspace*{1cm} \underline{\hspace{12cm}}

\[
\begin{array}{cc|cccc}
& & & $\hspace*{1cm} k$ & & \\
    &     &    $2$   &    $3$   &     $4$ & $5$ \\ \hline
    & $2$ & $0.7500$ & $0.8333$ & $0.8750$ & $0.9000$ \\
$r$ & $3$ & $0.5286$ & $0.6151$ & $0.6666$ & $0.7018$ \\
    & $4$ & $0.4047$ & $0.4800$ & $0.5275$ & $0.5614$ \\
    & $5$ & $0.3272$ & $0.3921$ &  $0.4343$ & $0.4650$ \\
\end{array}
\]
\begin{center}
\textbf{Table~4.} Approximate values of $g(k,\delta)$ in Theorem~\ref{thmB_new} for small $k$ and $\delta$.
\end{center}

\hspace*{1cm} \underline{\hspace{12cm}}

\end{document}